\documentclass[11pt]{article}
\usepackage{amscd,amsmath,amssymb}
\usepackage[mathscr]{eucal}

\makeatletter \@addtoreset{equation}{section}\makeatother

\title{\bf q-Analogues for Green functions for powers\qquad\qquad\qquad
of the invariant Laplacian in the unit disc}
\author{D. Shklyarov}
\date{}

\newtheorem{theorem}{Theorem}[section]
\newtheorem{lemma}[theorem]{Lemma}
\newtheorem{proposition}[theorem]{Proposition}
\newtheorem{corollary}[theorem]{Corollary}

\begin{document}
\large \maketitle

\makeatletter
\let\@thefnmark\relax
\@footnotetext{This research was supported in part by Award No UM1-2091 of
the US Civilian Research \& Development Foundation} \makeatother

\section{Introduction}

 In 1993 W.K.~Hayman and B.~Korenblum published among other results
explicit formulae for Green functions of powers of the Laplace operator in
the balls in ${\mathbb R}^n$ (see \cite{HK}). J. Peetre and
M. Engli\u{s} \cite{EP} have obtained analogous results for some
powers of the M\"obius-invariant Laplace operator in the unit ball in
${\mathbb C}^n$. In the particular case of the unit disc ${\mathbb U}
\subset {\mathbb C}$ they have presented explicit formulae for Green
functions of the powers $\Delta$, $\Delta^2$, $\Delta^3$, $\Delta^4$ of the
M\"obius-invariant (equi\-va\-lently, $SU(1,1)$-invariant) Laplace operator.
For
$\Delta$ and $\Delta^2$ the Green functions are
$$\frac{1}{4\pi} \ln t, $$
$$\frac{1}{16\pi}(-\ln t \ln(1-t)-2\mathrm{Li}_2(t)+\frac{\pi^2}{3}),$$
where $$t=\frac{(1-|z|^2)(1-|w|^2)}{|1-z\overline w|^2},$$ $z,w \in {\mathbb
U},$ $\mathrm{Li}_2(t)=\sum_{m=1}^{\infty} \frac{t^m}{m^2}$ is Euler's
dilogarithm.

 The aim of the present work is the computation of q-analogues for these
kernels. Namely, we will concern with the quantum unit disc which is a
homogeneous space of the quantum group $SU(1,1)$, and q is the parameter
used in the theory of quantum groups \cite{ChP}. Of course, all our formulae
become the classical ones at the limit $q\rightarrow1$.

 The statement of our main result (Theorem \ref{princresult}) is:
for any finite function $f$ in the quantum unit disc the following
equalities hold
\begin{equation}\label{delta1}
\Delta_q^{-1}f=\int \limits_{{\mathbb U}_q}{\mathbb G}_1fd\nu,
\end{equation}
\begin{equation}\label{delta2}
 \Delta_q^{-2}f=\int \limits_{{\mathbb U}_q}{\mathbb G}_2fd\nu,
\end{equation}
(here $\Delta_q$ and $\int \limits_{{\mathbb U}_q}\cdot d\nu$ are
q-analogues of the $SU(1,1)$-invariant Laplace operator and invariant
integral in the unit disc respectively, ${\mathbb G}_1$ and ${\mathbb
G}_2$ are certain kernels given by explicit formulae
(\ref{mathbbG1}),(\ref{mathbbG2}); precise definitions are to be found
below).

Let us outline ideas of the paper. In the classical case all the integral
operators we are interested in are intertwining (i.e., they commute with the
action of the group $SU(1,1)$ in spaces of functions in the unit disc). Thus
the kernels of these operators are functions in the simplest one:
\begin{equation}\label{yadro}
\frac{(1-z\bar{z})(1-\zeta\bar{\zeta})}
{(1-z\bar{\zeta})(1-\zeta\bar{z})}.
\end{equation}

In \cite{D3} an algebra was considered of kernels of integral operators in
the quantum disc, and q-analogues were obtained of integer negative powers
of (\ref{yadro}):
\begin{equation}\label{yadr}
G_{-l}=\{(1-\zeta\zeta^*)^{-l}(q^2z^*\zeta;q^2)_l(z
\zeta^*;q^2)_l(1-z^*z)^{-l}\},
\end{equation}
with $(a;q)_l\stackrel{\rm
def}=(1-a)\cdot(1-aq)\cdot\ldots\cdot(1-aq^{l-1})$. Furthermore, the kernels
$G_l$ may be defined for any $l\in{\mathbb C}$ by "analytic continuation" in
the parameter $q^{-2l}$ (see \cite{D3} or Section~2 of the present paper).

In the classical case "any" function of the kernel (\ref{yadro}) can be
expanded in integral by powers of this kernel using the Melline transform.
It turned out that in the quantum case as well as in the classical one the
kernels ${\mathbb G}_1$, ${\mathbb G}_2$ of the {\bf intertwining}
integral operators $\Delta_q^{-1}$, $\Delta_q^{-2}$ can be written
in the form
$$
{\mathbb G}_1=\int G_l d\sigma_1(l),
$$
$$
{\mathbb G}_2=\int G_l d\sigma_2(l),
$$ with some distributions
$d\sigma_1(l)$, $d\sigma_2(l)$. We find these integrals applying
the ope\-ra\-tors $\Delta_q^{-1}$, $\Delta_q^{-2}$ to a q-analogue
(denoted by $f_0$ in the text) of the delta-function at the centre
of the disc.

\bigskip

\section{\bf Statement of results}

First of all, we remind some notations and results on function theory in the
quantum unit disc (see \cite{D1}-\cite{D4}).

 Let $q \in (0;1)$. We impose the notation ${\rm Pol}({\mathbb C})_q$ for
the involutive algebra given by its generator $z$ and the unique commutation
relation
\begin{equation}\label{qr}
z^*z=q^2 zz^*+1-q^2.
\end{equation}

Let $y\stackrel{\rm def}=1-zz^*.$ It is straightforward that $y=y^*$,
\begin{equation}\label{qr1}
 z^*y=q^2 yz^*,\qquad  zy=q^{-2}yz,
\end{equation}
and any element $f \in {\rm Pol}({\mathbb C})_q$ admits a unique
decomposition
\begin{equation}\label{decomp}
f=\sum_{m>0} z^m\psi_m(y)+\psi_0(y)+\sum_{m>0} \psi_{-m}(y)z^{*m}.
\end{equation}

It is also not hard to show that the algebra ${\rm Pol}({\mathbb C})_q$
admits a unique up to unitary equivalence {\bf faithful} $*$-representation
by bounded operators in a Hilbert space and the spectrum of the operator,
corresponding to the element $y$, is the set $\{0\}\bigcup
q^{2\mathbb{Z}_+}$ (we shall use notation $y$ both for the element of the
polynomial algebra as well as for an indeterminate in the set
$q^{2\mathbb{Z}_+}$). This allows one to introduce the algebra $D({\mathbb
U})_q$ of finite functions in the quantum unit disc. By the definition it
consists of {\bf finite} series of the form (\ref{decomp}) with ${\rm supp }
\psi_m \subset q^{2\mathbb{Z}_+}$, \ ${\rm card} ({\rm supp }
\psi_m)<\infty$.

 The linear functional \cite[Theorem 3.5]{D2}
\begin{equation}\label{invint}
\int \limits_{{\mathbb U}_q}fd\nu\stackrel{\rm def}=(1-q^2)\sum_{m=o}^\infty
\psi_{0}(q^{2m})q^{-2m},
\end{equation}
where $f=\sum_{m>0} z^m\psi_m(y)+\psi_0(y)+\sum_{m>0} \psi_{-m}(y)z^{*m}\in
D({\mathbb U})_q$, is a q-analogue for the $SU(1,1)$-invariant integral in
the unit disc (see Section 3). We impose the notation $L^2(d\nu)_q$ for the
completion of $D({\mathbb U})_q$ with respect to the norm
\begin{equation}\label{norm}
 \|f\|\stackrel{\rm def}=\left({\int \limits_{{\mathbb
U}_q}f^*fd\nu}\right)^{\frac{1}{2}}.
\end{equation}

We need the well known first-order differential calculus over ${\rm
Pol}({\mathbb C})_q$. It is a ${\rm Pol}({\mathbb C})_q$-module
$\Omega^1({\mathbb C})_q$ given by its generators $dz$, $dz^*$, and the
commutation relations
\begin{equation}\label{fodc}
z\cdot dz=q^{-2}dz\cdot z, \ z^*\cdot dz^*=q^{2}dz^*\cdot z^*,
\ z^*\cdot dz=q^{2}dz\cdot z^*, \ z\cdot dz^*=q^{-2}dz^*\cdot z,
\end{equation}
and equipped with a linear map $d:{\rm Pol}({\mathbb C})_q
\rightarrow \Omega^1({\mathbb C})_q$ such that

\medskip

1. $ d:z\mapsto dz,\quad z^*\mapsto dz^*;$

\medskip

2. $ d(f_1f_2)=df_1\cdot f_2+f_1\cdot df_2$ for any $f_1,f_2\in{\rm
Pol}({\mathbb C})_q$ (the Leibniz rule).\medskip

The partial derivatives $\frac{\partial}{\partial z}$,
$\frac{\partial}{\partial z^*}$ in ${\rm Pol}({\mathbb C})_q$ are given by
$$ df=dz\frac{\partial f}{\partial z}+dz^*\frac{\partial f}{\partial z^*}.$$

The operator
\begin{equation}\label{lapl}
\Delta_q f\stackrel{\rm def}{=}(1-zz^*)^2 \frac{\partial }{\partial
z^*}\frac{\partial }{\partial z}f
\end{equation}
is a q-analogue of the invariant Laplacian in the unit disc (see
\cite{D2,D4}).

 One can define the operators $\frac{\partial}{\partial z}$,
$\frac{\partial}{\partial z^*}$, $\Delta_q$ on the space of finite functions
(it is sufficient to use the formulae
$$df(y)=-q^2\frac{f(y)-f(q^2y)}{y-yq^2}z^*dz-z\frac{f(y)-f(q^2y)}{y-yq^2}
dz^*,$$
$$ f(y)dz=dzf(y),$$ $$ f(y)dz^*=dz^*f(y),$$ for $f\in {\rm Pol}({\mathbb
C})_q$).

The following result was announced in \cite[Proposition 3.2]{D1}:\medskip

\begin{theorem}\label{invertible}
The operator $\Delta_q$ can be extended to the selfadjoint bounded
invertible operator in $L^2(d\nu)_q$.
\end{theorem}

To formulate our results we need the notion of an integral operator in the
quantum case.

Impose the notation $D({\mathbb U}\times{\mathbb U})'_q$ for the space of
formal series of the form
\begin{equation}\label{kernel}
f=\sum_{(i,j)\in{\mathbb Z}^2}f_{ij},
\end{equation}
\begin{equation}\label{kernell}f_{ij}=
\left\{
\begin{array}{lr}
z^i\otimes1 \cdot \psi_{ij}(y\otimes1,1\otimes y)\cdot 1\otimes z^j, &
\mbox{$i\geq0,j\geq0, $}\\
z^i\otimes1 \cdot \psi_{ij}(y\otimes1,1\otimes y)\cdot 1\otimes z^{*j}, &
\mbox{$i\geq0,j<0,  $}\\
z^{*i}\otimes1\cdot \psi_{ij}(y\otimes1,1\otimes y)\cdot 1\otimes z^j,  &
\mbox{$i<0,j\geq0,  $}\\
z^{*i}\otimes1\cdot \psi_{ij}(y\otimes1,1\otimes y)\cdot 1\otimes z^{*j},&
\mbox{$i<0,j<0,  $}
\end{array}
\right.
\end{equation}
with $\{\psi_{ij}\}$ being any functions on $q^{2{\mathbb Z}_+}\times
q^{2{\mathbb Z}_+}$. It is convenient in the sequel to write
$z,y,z^*,\zeta,\eta,\zeta^*$ instead of
$z\otimes1,y\otimes1,z^*\otimes1,1\otimes z,1\otimes y,1\otimes
z^*$, respectively. Note that $D({\mathbb U}\times{\mathbb U})'_q$
can be made into a topological vector space. We describe the
topology in Section 3.

Let $D({\mathbb U})'_q$ be the space of formal series of the form
(\ref{decomp}) with $\{\psi_m(y)\}_{m\in{\mathbb Z}}$ being any functions on
$q^{2{\mathbb Z}_+}$. Then one shows that for any $K\in D({\mathbb
U}\times{\mathbb U})'_q$ the map
$$\int \limits_{{\mathbb
U}_q}Kfd\nu\stackrel{\rm def}=id\otimes\nu\left(K\cdot1\otimes
f\right),\quad f\in D({\mathbb U})_q,\quad \nu(f)\stackrel{\rm def}=\int
\limits_{{\mathbb U}_q}fd\nu\,,
$$ is a well defined operator from $D({\mathbb
U})_q$ into $D({\mathbb U})'_q$. (Indeed, it follows from the relations
$$z^*\phi(y)=\phi(q^2y)z^*,\qquad z\phi(y)=\phi(q^{-2}y)z,\qquad {\rm supp}
\phi\subset q^{2{\mathbb Z}_+}$$
and
$$z^{*k}z^k=(1-q^2y)(1-q^4y)\dots(1-q^{2k}y), \ z^kz^{*k}=(1-y)(1-q^{-2}y)
\dots(1-q^{-2k+2}y),$$
that $D({\mathbb U})'_q$ is a $D({\mathbb U})_q$-bimodule. Now to prove the
correctness of the definition it sufficies to observe that for a finite
$f_2(y)$ $$\int \limits_{{\mathbb U}_q}z^kf_1(y)z^jf_2(y)d\nu=0,\quad k\ne0
\quad {\rm or}\quad j\ne0,$$ $$\int \limits_{{\mathbb
U}_q}f_1(y)z^{*k}z^jf_2(y)d\nu=0,\quad k\ne j,$$
$$\int \limits_{{\mathbb U}_q}f_1(y)z^{*k}f_2(y)z^{*j}d\nu=0,\quad k\ne0
\quad {\rm
or}\quad j\ne0,$$ $$\int \limits_{{\mathbb
U}_q}z^kf_1(y)f_2(y)z^{*j}d\nu=0,\quad k\ne j.)$$

Such operators can be treated as integral operators.

We have already mentioned (see (\ref{yadr})) one important set of kernels
introduced in \cite[Section 6]{D3}:

\begin{equation}\label{yadr1}
G_{-l}=\{(1-\zeta\zeta^*)^{-l}(q^2z^*\zeta;q^2)_l(z
\zeta^*;q^2)_l(1-z^*z)^{-l}\},\qquad l=1,2,3\ldots\,,
\end{equation}
where $(a;q)_n\stackrel{\rm def}=(1-a)(1-aq)\dots(1-aq^{n-1})$, and the
brackets $\{\}$ indicate that $z,z^*$ should be multiplied within the
algebra ${\rm Pol}({\mathbb C})^{op}_q$ derived from ${\rm Pol}({\mathbb
C})_q$ by replacing its product to the opposite one (for example, $\{z\cdot
z^*\}=z^*\cdot z$). These kernels are invariant in a sense (we explain the
term "invariant" in Section 3; see also \cite{D3}) and therefore may be
regarded as q-analogues of the kernels
\begin{equation}\label{clasyadr}
\left(\frac{(1-z\bar{z})(1-\zeta\bar{\zeta})}
{(1-z\bar{\zeta})(1-\zeta\bar{z})}\right)^{-l},
\end{equation}
with $l$ being a positive integer number. Using relation (1.3.2) from
\cite{GR} one can rewrite
\begin{equation}\label{ker1}
G_{-l}=\sum_{k=0}^{\infty}\sum_{n=0}^{\infty}q^{2k}\frac{(q^{-2l};q^{2})_k
\cdot(q^{-2l};q^{2})_n}{(q^{2};q^{2})_k\cdot(q^{2};q^{2})_n}z^ny^{-l}z^{*k}
\zeta^k\eta^{-l}\zeta^{*n},
\end{equation}
(the sum is finite). It is evident that for any finite functions
$\phi_1,\phi_2$ the function $t\stackrel{\rm def}=q^{-2l}\mapsto
g_{\phi_1,\phi_2}(t)\stackrel{\rm def}=\int \limits_{{\mathbb
U}_q}\phi_1\cdot\left(\int \limits_{{\mathbb
U}_q}G_{-l}\cdot\phi_2d\nu\right)d\nu$ belongs to ${\mathbb C}[t,t^{-1}]$.
This observation allows one to prove (just as it was done in \cite{D3}) that
there exists a unique vector-function of a {\bf complex} variable $t$ with
values in $D({\mathbb U}\times{\mathbb U})'_q$ which coincides with the
right side of (\ref{ker1}) for $t\in q^{-2\mathbb N}$. In the sequel $G_l$
will stand for this "analytic continuation". In this way one obtains
q-analogues of (\ref{clasyadr}) for any complex power $l$.

We will need also the kernels
\begin{equation}\label{ker2}
\hat{G}_N\stackrel{\rm def}=\lim_{l\rightarrow
N}\frac{G_l-G_N}{l-N},\qquad N=1,2,\ldots
\end{equation}
(the limit in the topology mentioned above).\smallskip

{\bf Remark.} Let
$$
L_0(\xi)=0,\quad L_k(\xi)\stackrel{\rm
def}=\frac{1}{1-\xi}+\frac{q^2}{1-q^2\xi}+\dots+\frac{q^{2k-2}}{1-q^{2k-2}
\xi}, \quad k=1,2, \ldots , \infty.
$$
Using (\ref{ker1}) and
formulae
$$\frac{d}{dl}y^l=\ln y\cdot y^l,\quad
\frac{d}{dl}(q^{2l};q^2)_n=h\cdot L_n(q^{2l})\cdot
(q^{2l};q^2)_n\,,
$$ one can show that for $h\stackrel{\rm def}=\ln
q^{-2}$
\begin{equation}\label{ker3}
\hat{G}_N=\sum_{k=0}^{\infty}\sum_{n=0}^{\infty}q^{2k}\frac{(q^{2N};q^{2})_k
\cdot(q^{2N};q^{2})_n}{(q^{2};q^{2})_k\cdot(q^{2};q^{2})_n}\cdot
\Psi_{N,k,n},
\end{equation}
where
$$\Psi_{N,k,n}\stackrel{\rm def}=h(L_k(q^{2N})+L_n(q^{2N}))z^ny^Nz^{*k}
\zeta^k\eta^N\zeta^{*n}+ z^ny^N\ln(y)z^{*k}\zeta^k\eta^N\zeta^{*n}
$$
$$+z^ny^Nz^{*k}\zeta^k\eta^N\ln(\eta)\zeta^{*n}.$$

Note that
$$\lim_{q\rightarrow1}G_l=\left(\frac{(1-z\bar{z})(1-\zeta\bar{\zeta})}
{(1-z\bar{\zeta})(1-\zeta\bar{z})}\right)^l,$$
$$\lim_{q\rightarrow1}\hat{G}_N=\left(\frac{(1-z\bar{z})(1-\zeta\bar
{\zeta})}
{(1-z\bar{\zeta})(1-\zeta\bar{z})}\right)^N\cdot
\ln\left(\frac{(1-z\bar{z})(1-\zeta\bar{\zeta})}
{(1-z\bar{\zeta})(1-\zeta\bar{z})}\right).$$

The principal result of the present work is

\begin{theorem}\label{princresult}
For any $f\in D({\mathbb U})_q$
\begin{equation}\label{delta^-1}
\Delta_q^{-1}f=\int \limits_{{\mathbb U}_q}{\mathbb G}_1fd\nu,
\end{equation}
\begin{equation}\label{delta^-2}
 \Delta_q^{-2}f=\int \limits_{{\mathbb U}_q}{\mathbb G}_2fd\nu,
\end{equation}
where
\begin{equation}\label{mathbbG1}
 {\mathbb G}_1\stackrel{\rm
def}=-\sum_{m=1}^{\infty}\frac{q^{-2}-1}{q^{-2m}-1}G_m\,,
\end{equation}
\begin{equation}\label{mathbbG2}
 {\mathbb G}_2\stackrel{\rm
def}=\sum_{m=1}^{\infty}\frac{q^{2m-2}(1+q^{2m})(1-q^2)^2}{(1-q^{2m})^2}G_m-
\frac{1-q^2}{h}\sum_{m=1}^{\infty}\frac{q^{-2}-1}{q^{-2m}-1}\hat{G}_m\,,
\end{equation}
and $G_m,\hat{G}_m$ are given by (\ref{ker1}), (\ref{ker2}), respectively.
\end{theorem}

\section{\bf Auxiliary result: radial part of the invariant Laplacian}

It is easy to check by direct calculations that for any $f(y)\in{\rm
Pol}({\mathbb C})_q$ or $D({\mathbb U})_q$ $$\Delta_q f(y)=q^{-1}y^2
D(1-qy)Df(y),$$ where $$(Df)(t)\stackrel{\rm
def}=\frac{f(q^{-1}t)-f(qt)}{q^{-1}t-qt}.$$

Let $L^2(d\nu)^{(0)}_q\stackrel{\rm def}=\{f(y)\in D({\mathbb
U})'_q|\sum_{m=0}^{\infty}|f(q^{2m})|^2q^{-2m}<\infty\}.$ The following
proposition is proved in \cite[Lemma 5.5]{D2}:

\begin{proposition}\label{delta0}
$\Delta_q^{(0)}\stackrel{\rm def}=q^{-1}y^2 D(1-qy)D$ is a bounded
selfadjoint invertible operator in $L^2(d\nu)^{(0)}_q$.
\end{proposition}

The term "radial part of the invariant Laplacian in the quantum disc" stand
for this operator.

Let $f_0=f_0(y)$ be such a finite function that
\begin{equation}\label{f-0}
f_0=
  \left\{
\begin{array}{ll}
1 & \mbox{$y=1,$}\\
0 & \mbox{$y=q^{2k},k=1,2, \ldots . $}
\end{array}
\right.
\end{equation}

In this section we will prove the following

\begin{theorem}\label{g1,g2}
$$(\Delta_q^{(0)})^{-1}f_0=g_1(y),$$ $$(\Delta_q^{(0)})^{-2}f_0=g_2(y),$$
where
\begin{equation}\label{g1}
 g_1(y)=-(1-q^2)\sum_{m=1}^{\infty}\frac{q^{-2}-1}{q^{-2m}-1}y^m,
\end{equation}
$$
g_2(y)=(1-q^2)
$$
\begin{equation}\label{g2}
\times\left(\sum_{m=1}^{\infty}\frac{q^{2m-2}(1+q^{2m})(1-q^2)^2}
{(1-q^{2m})^2}y^m- \frac{1-q^2}{h}\ln y
\sum_{m=1}^{\infty}\frac{q^{-2}-1}{q^{-2m}-1}y^m\right).
\end{equation}
\end{theorem}

Remind some well known notations \cite{GR}:
$$(a;q)_n\stackrel{\rm def}=(1-a)(1-aq)\dots(1-aq^{n-1}),$$
$$(a;q)_{\infty}\stackrel{\rm
def}=(1-a)(1-aq)\dots(1-aq^{n-1})\dots\,,$$
$$\Gamma_q(x)=\frac{(q;q)_\infty}{(q^x;q)_\infty}(1-q)^{1-x},$$
$$_r \Phi_s \left[{a_1,a_2,\ldots,a_r;q;z \atop b_1,\ldots,b_s}\right]$$
$$=\sum_{n \in{\mathbb Z}_+} \frac{(a_1;q)_n \cdot(a_2;q)_n \cdot \ldots
\cdot(a_r;q)_n} {(b_1;q)_n \cdot(b_2;q)_n \cdot \ldots
\cdot(b_s;q)_n(q;q)_n}\left((-1)^n \cdot
q^\frac{n(n-1)}{2}\right)^{1+s-r}\cdot z^n,$$
$$\int\limits_{0}^{1}f(y)d_{q^2}y=(1-q^2)\sum_{m=0}^{\infty}f(q^{2m})q^{2m}.
$$

We will use the following results from \cite{D1,D2,D4}:

\begin{proposition}\label{spher}
The functions
\begin{equation}\label{sph}
 \varphi_{\rho}(y)\stackrel{\rm def}=_3\Phi_2 \left[{y^{-1},q^{1+2i\rho},
q^{1-2i\rho};q^2;q^2 \atop q^2,0}\right],
\end{equation}
\begin{equation}\label{jost}
 \psi_{\rho}(y)\stackrel{\rm def}=
y^{\frac{1}{2}-i\rho}\cdot _2\Phi_1
\left[{q^{1-2i\rho},q^{1-2i\rho};q^2;q^2y \atop q^{2-4i\rho}}\right],
\end{equation}
and $\psi_{-\rho}(y)$ for $\rho\in{\mathbb C} \backslash\frac{1}{2i}{\mathbb
N}$ are solutions of the equation
$$\Delta_q^{(0)}f(y)=\lambda(\rho)f(y),$$ where
$$\lambda(\rho)=-\frac{(1-q^{1+2i\rho})(1-q^{1-2i\rho})}{(1-q^2)^2}.$$
Moreover,
\begin{equation}\label{linearcomb}
\varphi_\rho(y)=\frac{\Gamma_{q^2}(2i\rho)}{\Gamma^2_{q^2}(\frac{1}{2}+i
\rho)}
\psi_\rho(y)+\frac{\Gamma_{q^2}(-2i\rho)}{\Gamma^2_{q^2}(\frac{1}{2}-i\rho)}
\psi_{-\rho}(y).
\end{equation}
\end{proposition}
\smallskip

{\bf Remark.} $\varphi_{\rho}(y)$ is a q-analogue of the spherical function
in the unit disc (see \cite{Helg}).\medskip

\begin{proposition}\label{spectr}\cite[Corollary 4.2]{D4}.
The spectrum of
$\Delta_q^{(0)}$ is simple purely continuous and coincides with the segment
$$\left[-\frac{1}{(1-q)^2};-\frac{1}{(1+q)^2}\right].$$
\end{proposition}
\medskip

\begin{proposition}  \label{borel}\cite[Proposition 4.17]{D4}.
Consider the Borel
measure $d\sigma$ on the segment $[0;\frac{2\pi}{h}]$ ($h=\ln q^{-2}$) given
by
\begin{equation}\label{dsigma}
 d\sigma(\rho)=\frac{1}{4\pi}\frac{h}{1-q^2}\frac{\Gamma^2_{q^2}
(\frac{1}{2}-i\rho)
\Gamma^2_{q^2}(\frac{1}{2}+i\rho)}{\Gamma_{q^2}(-2i\rho)\Gamma_{q^2}
(2i\rho)}d\rho.
\end{equation}
Then the linear operator
\begin{equation}\label{sphtrans}
 f(y)\mapsto\hat{f}(\rho)\stackrel{\rm
def}=\int\limits_{0}^{1}\varphi_\rho(y)f(y)y^{-2}d_{q^2}y,
\end{equation}
defined on functions with finite supports inside $q^{2{\mathbb Z}_+}$ is
extendable by continuity to a unitary operator
$$u:L^2(d\nu)^{(0)}_q\rightarrow L^2(d\sigma).$$
For all $f\in L^2(d\nu)^{(0)}_q$
\begin{equation}\label{mult by fun}
 u\cdot\Delta_q^{(0)}f=\lambda(\rho)uf,
\end{equation}
and the inverse operator is
\begin{equation}\label{inv sph tr}
 f(\rho)\mapsto \int\limits_{0}^{2\pi/h}\varphi_\rho(y)f(\rho)d\sigma(\rho).
\end{equation}
\end{proposition}\medskip

{\bf Remark.} Formulae (\ref{sphtrans}), (\ref{inv sph tr}) present a
decomposing in eigenfunctions of the operator $\Delta_q^{(0)}$. The function
$\hat{f}(\rho)$ is called the spherical transform of $f(y)$ while $f(y)$ is
the inverse spherical transform for $\hat{f}(\rho)$.

Now let us turn to proving of Theorem \ref{g1,g2}

\begin{lemma}\label{gmq2n} Let $g_m(y)$ stand for the function such that
\begin{equation}\label{defgm}
\Delta_q^{(0)m}g_m(y)=f_0.
\end{equation}
Then
$$ g_m(q^{2N})=
(-1)^m(1-q^2)^{2m}(q^2;q^2)_{\infty}q^N\sum_{k=0}^{\infty} \frac{q^{2Nk+2k}}
{(q^2;q^2)_k}$$
\begin{equation}\label{explgm}
\times{\rm Res}_{\tau=q}\left(\frac{\tau^{N+m-1}(1-\tau^2)
(q^{2k+2}\tau^2;q^2)_\infty
d\tau}{(\tau-q)^m(1-q\tau)^m(q^{2k+1}\tau;q^2)^2_ \infty}\right).
\end{equation}
\end{lemma}

{\bf Proof of the lemma.} Applying the spherical transform to the both sides
of (\ref {defgm}), using (\ref {mult by fun}) and equality
$\varphi_\rho(1)=1$ we get
\begin{equation}\label{hatgm} \lambda(\rho)^m\hat{g}_m(\rho)=1-q^2
\end{equation}
and then
\begin{equation}\label{hatgm1}
 \hat{g}_m(\rho)=(-1)^m\frac{(1-q^2)^{2m+1}}{(1-q^{1+2i\rho})^m
(1-q^{1-2i\rho})^m}\,.
\end{equation}

Now to obtain $g_m(y)$ in an explicit form it is sufficient to apply the
inverse spherical transform to the both sides of (\ref {hatgm1}), i.e.,
\begin{equation}\label{gm(y)}
 g_m(y)=(-1)^m(1-q^2)^{2m+1}\int\limits_{0}^{2\pi/h}
\frac{\varphi_\rho(y)}{(1-q^{1+2i\rho})^m(1-q^{1-2i\rho})^m}d\sigma(\rho).
\end{equation}

Next, to compute the integral in the right side of (\ref{gm(y)}) we replace
$\varphi_\rho(y)$ by its decomposition into sum of two items
(cf.(\ref{linearcomb})):
$$g_m(y)=c_m\int\limits_{0}^{2\pi/h}
\frac{\psi_ \rho(y)}{(1-q^{1+2i\rho})^m(1-q^{1-2i\rho})^m}
\frac{\Gamma^2_{q^2}(\frac{1}{2}-i\rho)}{\Gamma_{q^2}(-2i\rho)}d\rho $$
\begin{equation}\label{twoint}
+ c_m\int\limits_{0}^{2\pi/h}
\frac{\psi_{-\rho}(y)}{(1-q^{1+2i\rho})^m(1-q^{1-2i\rho})^m}
\frac{\Gamma^2_{q^2}(\frac{1}{2}+i\rho)}{\Gamma_{q^2}(2i\rho)}d\rho,
\end{equation}
where $c_m=\frac{(-1)^m}{4\pi}(1-q^2)^{2m}h$.

 The two integrals in the right side of (\ref{twoint}) are equal to
each other (to check this one should replace $\rho$ by $-\rho$ in the former
integral and observe that all the functions under the integrals are
$\frac{2\pi}{h}$-periodic).

Hence
\begin{equation}\label{oneterm}
g_m(y)=\frac{(-1)^m}{2\pi}(1-q^2)^{2m}h\int\limits_{0}^{2\pi/h}
\frac{\psi_\rho(y)}{(1-q^{1+2i\rho})^m(1-q^{1-2i\rho})^m}
\frac{\Gamma^2_{q^2}(\frac{1}{2}-i\rho)}{\Gamma_{q^2}(-2i\rho)}d\rho\,.
\end{equation}

Now let us make use of the equalities
\begin{equation}\label{psirho}
 \psi_\rho(y)=y^{1/2-i\rho}\sum_{k=0}^{\infty}\frac{(q^{1-2i\rho};q^2)^
{2}_k}
{(q^{2-4i\rho};q^2)_k(q^2;q^2)_k}q^{2k}y^k
\end{equation}
and
\begin{equation}\label{gam/gam}
\frac{\Gamma^2_{q^2}(\frac{1}{2}-i\rho)}{\Gamma_{q^2}(-2i\rho)}=
\frac{\frac{(q^2;q^2)^2_\infty}{(q^{1-2i\rho};q^2)^2_\infty}(1-q^2)^
{1+2i\rho}}
{\frac{(q^2;q^2)_\infty}{(q^{-4i\rho};q^2)_\infty}(1-q^2)^{1+2i\rho}}=
(q^2;q^2)_\infty\frac{(q^{-4i\rho};q^2)_\infty}{(q^{1-2i\rho};q^2)^2_
\infty}\,.
\end{equation}
Thus $$g_m(y)=\frac{(-1)^mh}{2\pi}(1-q^2)^{2m}(q^2;q^2)_\infty
y^{1/2}\sum_{k=0}^{\infty}\frac{q^{2k}y^k}{(q^2;q^2)_k}$$
$$\times\int\limits_{0}^{2\pi/h}
\frac{1}{(1-q^{1+2i\rho})^m(1-q^{1-2i\rho})^m}\frac{(q^{1-2i\rho};q^2)^2_k}
{(q^{2-4i\rho};q^2)_k}\frac{(q^{-4i\rho};q^2)_\infty}{(q^{1-2i\rho};q^2)^
2_\infty}
y^{-i\rho}d\rho$$ (we have exchanged summation over $k$ and integration over
$[0;\frac{2\pi}{h}]$ because of the uniform convergence of the series
(\ref{psirho}) for any fixed $y\in q^{2{\mathbb Z}_+}$).

Let $y=q^{2N}$. Remind that $q=e^{-h/2}$. Then
$$\int\limits_{0}^{2\pi/h}
\frac{1}{(1-q^{1+2i\rho})^m(1-q^{1-2i\rho})^m}\frac{(q^{1-2i\rho};q^2)^2_k}
{(q^{2-4i\rho};q^2)_k}\frac{(q^{-4i\rho};q^2)_\infty}{(q^{1-2i\rho};q^2)^
2_\infty}
q^{-2Ni\rho}d\rho$$
$$=\int\limits_{0}^{2\pi/h}
\frac{1}{(1-q e^{-hi\rho})^m(1-q e^{hi\rho})^m}\frac{(q e^{hi\rho};q^2)^2_k}
{(q^2 e^{2hi\rho};q^2)_k}\frac{(e^{2hi\rho};q^2)_\infty}{(q
e^{hi\rho};q^2)^2_\infty} e^{hNi\rho}d\rho$$
$$=\frac{1}{hi}\int\limits_{\mathbb T}\frac{1}{(1-\frac{q}{\tau})^m(1-q\tau)
^m}
\frac{(q\tau;q^2)^2_k} {(q^2
\tau^2;q^2)_k}\frac{(\tau^2;q^2)_\infty}{(q\tau;q^2)^2_\infty}
\tau^{N-1}d\tau$$
$$=\frac{1}{hi}\int\limits_{\mathbb T}\frac{\tau^{N+m-1}}{(\tau-q)^m
(1-q\tau)^m}
(1-\tau^2)\frac{(q^{2k+2}\tau^2;q^2)_\infty}{(q^{2k+1} \tau;q^2)^2_\infty}
d\tau\,.$$ This completes the proof of the lemma.\hfill $\square$
\medskip

Thus by Lemma \ref{gmq2n}
 we have $$g_1(q^{2N})=
-(1-q^2)^{2}(q^2;q^2)_{\infty}q^{2N}\sum_{k=0}^{\infty} \frac{q^{2Nk+2k}}
{(q^2;q^2)_k}\cdot\frac{(q^{2k+4};q^2)_\infty}{(q^{2k+2};q^2)^2_\infty}$$
$$=-(1-q^2)^{2}(q^2;q^2)_{\infty}q^{2N}\sum_{k=0}^{\infty} \frac{q^{2Nk+2k}}
{(q^2;q^2)_k(1-q^{2k+2})(q^{2k+2};q^2)_\infty}$$
$$=-(1-q^2)^{2}(q^2;q^2)_{\infty}q^{2N}\sum_{k=0}^{\infty} \frac{q^{2Nk+2k}}
{(1-q^{2k+2})(q^{2};q^2)_\infty}.$$ Hence
\begin{equation}\label{g1y}
 g_1(y)=-(1-q^2)\sum_{m=1}^{\infty}\frac{q^{-2}-1}{q^{-2m}-1}y^m.
\end{equation}

For $g_2$ the calculations are much more complicated
$${\rm Res}_{\tau=q}\left(\frac{\tau^{N+1}(1-\tau^2)
(q^{2k+2}\tau^2;q^2)_\infty
d\tau}{(\tau-q)^2(1-q\tau)^2(q^{2k+1}\tau;q^2)^2_ \infty}\right)$$
$$=\frac{d}{d\tau}\left(\frac{\tau^{N+1}(1-\tau^2)}{(1-q\tau)^2}
\cdot(q^{2k+2}\tau^2;q^2)_\infty\frac{1}{(q^{2k+1}\tau;q^2)^2_
\infty}\right)_{\tau=q}.$$
One proves that
$$\frac{d}{d\tau}(\tau;q^2)_\infty=-(\tau;q^2)_\infty\cdot L_\infty(\tau)$$
(here we use the notation $L_\infty$ from Section 2). Thus
$$\frac{d}{d\tau}\left(\frac{\tau^{N+1}(1-\tau^2)}{(1-q\tau)^2}
\cdot(q^{2k+2}\tau^2;q^2)_\infty\frac{1}{(q^{2k+1}\tau;q^2)^2_
\infty}\right)$$
$$=\frac{(N+1)\tau^N-(N+3)\tau^{N+2}}{(1-q\tau)^2}\cdot\frac
{(q^{2k+2}\tau^2;q^2)_\infty}
{(q^{2k+1}\tau;q^2)^2_\infty}$$
$$+\frac{2q(\tau^{N+1}-\tau^{N+3})}
{(1-q\tau)^3}\cdot\frac{(q^{2k+2}\tau^2;q^2)_\infty}
{(q^{2k+1}\tau;q^2)^2_\infty}+\frac{\tau^{N+1}-\tau^{N+3}} {(1-q\tau)^2}$$
$$\times\left(-2q^{2k+2}\tau
L_\infty(q^{2k+2}\tau^2)\frac {(q^{2k+2}\tau^2;q^2)_\infty}
{(q^{2k+1}\tau;q^2)^2_\infty}+ 2q^{2k+1}L_\infty(q^{2k+1}\tau)\frac
{(q^{2k+2}\tau^2;q^2)_\infty}
{(q^{2k+1}\tau;q^2)^2_\infty}\right).$$
For $\tau=q$ we get
$$\frac{(N+1)q^N-(N+3)q^{N+2}}{(1-q^2)^2}\cdot\frac
{(q^{2k+4};q^2)_\infty} {(q^{2k+2};q^2)^2_\infty}+\frac{2q(q^{N+1}-q^{N+3})}
{(1-q^2)^3}\cdot\frac{(q^{2k+4};q^2)_\infty} {(q^{2k+2};q^2)^2_\infty}$$
$$-\frac{q^{N+1}-q^{N+3}}
{(1-q^2)^2}\cdot2q^{2k+3} L_\infty(q^{2k+4})\frac {(q^{2k+4};q^2)_\infty}
{(q^{2k+2};q^2)^2_\infty}$$
$$+\frac{q^{N+1}-q^{N+3}}{(1-q^2)^2}
2q^{2k+1}L_\infty(q^{2k+2})\frac {(q^{2k+4};q^2)_\infty}
{(q^{2k+2};q^2)^2_\infty}$$
$$=q^N\frac{(N+1)-(N+3)q^{2}}{(1-q^2)^2}\cdot\frac
{1} {(1-q^{2k+2})(q^{2k+2};q^2)_\infty}$$
$$+\frac{2q^{N+2}}
{(1-q^2)^2}\cdot\frac{1} {(1-q^{2k+2})(q^{2k+2};q^2)_\infty}$$
$$-\frac{q^{N+1}}
{(1-q^2)}\cdot2q^{2k+3} L_\infty(q^{2k+4})\frac {1}
{(1-q^{2k+2})(q^{2k+2};q^2)_\infty}$$ $$+\frac{q^{N+1}}{(1-q^2)}
2q^{2k+1}L_\infty(q^{2k+2})\frac {1} {(1-q^{2k+2})(q^{2k+2};q^2)_\infty}$$
$$=\frac {1} {(1-q^2)(1-q^{2k+2})(q^{2k+2};q^2)_\infty}$$
$$\times\left(q^N(N+1)
+2q^{2k+1}\cdot q^{N+1}(L_\infty(q^{2k+2})-q^2L_\infty(q^{2k+4}))\right).$$
But
$$L_\infty(q^{2k+2})-q^2L_\infty(q^{2k+4})=\sum_{m=0}^{\infty}\frac{q^{2m}}
{1-q^{2m+2k+2}}-\sum_{m=0}^{\infty}\frac{q^{2m+2}}
{1-q^{2m+2k+4}}=\frac{1}{1-q^{2k+2}}\,.$$
Hence finally we obtain
$${\rm Res}_{\tau=q}\left(\frac{\tau^{N+1}(1-\tau^2)
(q^{2k+2}\tau^2;q^2)_\infty
d\tau}{(\tau-q)^2(1-q\tau)^2(q^{2k+1}\tau;q^2)^2_ \infty}\right)$$
$$=\frac {1}
{(1-q^2)(1-q^{2k+2})^2(q^{2k+2};q^2)_\infty}\cdot(q^N(N+1)(1-q^{2k-2})
+2q^N\cdot q^{2k+2})\,,$$ and, using (\ref{explgm}),
$$ g_2(q^{2N})=
(1-q^2)^{3}(q^2;q^2)_{\infty}q^{2N}\sum_{k=0}^{\infty} \frac{q^{2Nk+2k}}
{(q^2;q^2)_k(1-q^{2k+2})^2(q^{2k+2};q^2)_\infty}$$
$$\times((N+1)(1-q^{2k+2})+2q^{2k+2}).$$
Now to complete the proof of Theorem \ref{g1,g2}
 it is sufficient to replace $q^{2N}$ by $y$ and $N$ by
$-\frac{1}{h}\ln y$ in the last formula.\smallskip

{\bf Remark. }

1. The explicit form (\ref{g1}) of $g_1(y)$ was obtained in
\cite[Proposition 1.1]{D4} in another way.

2. The function $g_2(y)$ given by (\ref{g2}) can be treated as a q-analogue
of Rogers' dilogarithm.

\section{\bf Some more auxiliary results: quantum symmetry}

Remind that the quantum universal enveloping algebra $U_q \mathfrak{sl}_2$
is a Hopf algebra over ${\mathbb C}$ determined by the generators $K,K^{-1},
E, F$ and the relations
$$KK^{-1}=K^{-1}K=1,\quad K^{\pm 1}E=q^{\pm 2}EK^{\pm 1},\quad K^{\pm
1}F=q^{\mp 2}FK^{\pm 1},$$
$$EF-FE=(K-K^{-1})/(q-q^{-1}),$$
$$\Delta(K^{\pm 1})=K^{\pm 1}\otimes K^{\pm 1},\quad \Delta(E)=E \otimes 1+K
\otimes E,\quad \Delta(F)=F \otimes K^{-1}+1 \otimes F.$$

 Note that
$$\varepsilon(E)=\varepsilon(F)=\varepsilon(K^{\pm 1}-1)=0,$$
$$S(K^{\pm 1})=K^{\mp 1},\quad S(E)=-K^{-1}E,\quad S(F)=-FK,$$
with $\varepsilon:U_q \mathfrak{sl}_2 \to{\mathbb C}$ and $S:U_q
\mathfrak{sl}_2 \to U_q \mathfrak{sl}_2$ being respectively the counit and
the antipode of $U_q \mathfrak{sl}_2$.

Let $F$ stand for an algebra over ${\mathbb C}$ with a unit and equipped
also with a structure of $U_q \mathfrak{sl}_2$-module. $F$ is called an $U_q
\mathfrak{sl}_2$-module (covariant) algebra if

1. the multiplication $m:F\otimes F\rightarrow F$ is a morphism of $U_q
\mathfrak{sl}_2$-modules;

2. for any $\xi\in U_q \mathfrak{sl}_2$
$$\xi(1)=\varepsilon(\xi)\cdot1$$ (here $1$ is the unit
of $F$). Note that an element $v$ of an $U_q \mathfrak{sl}_2$-module is
called invariant if for any $\xi\in U_q \mathfrak{sl}_2$
$$\xi(v)=\varepsilon(\xi)\cdot v.$$

Let $M$ be an $U_q \mathfrak{sl}_2$-module and $F$-bimodule for some
covariant algebra $F$. Then $M$ is called covariant if the multiplication
maps
$$m_L:F\otimes M\rightarrow M,\quad m_R:M\otimes F\rightarrow M$$
are morphisms of $U_q \mathfrak{sl}_2$-modules.

Equip $U_q \mathfrak{sl}_2$ with the involution given by
\begin{equation}\label{invol}
 E^*=-KF,\quad F^*=-EK^{-1},\quad (K^{\pm 1})^*=K^{\pm 1}.
\end{equation}
$U_q \mathfrak{su}_{1,1}$ is the $*$-Hopf algebra produced this way.

An involutive algebra $F$ is said to be $U_q \mathfrak{su}_{1,1}$-module
algebra (covariant $*$-algebra) if it is an $U_q \mathfrak{sl}_{2}$-module
one and
$$(\xi f)^*=(S(\xi))^*\cdot f^*$$
for any $\xi\in U_q \mathfrak{su}_{1,1}$ and $f\in F$.

It is very well known (see, for instance, \cite{D2}) that ${\rm Pol}(
{\mathbb
C})_q$ can be equipped with a structure of a covariant $*$-algebra in the
following way:
\begin{equation}\label{action1}
 K^{\pm1}z=q^{\pm2}z,\quad Ez=-q^{1/2}z^2,\quad Fz=q^{1/2},
\end{equation}
\begin{equation}\label{action2}
 K^{\pm1}z^*=q^{\mp2}z^*,\quad Ez^*=q^{-3/2},\quad Fz=-q^{-5/2}z^{*2}.
\end{equation}

The formulae (\ref{action1}),(\ref{action2}) imply: for any polynomial $f$
\begin{equation}\label{action3}
 K^{\pm1}f(y)=f(y),\quad Ef(y)=-\frac{q^{1/2}}{1-q^2}z(f(y)-f(q^2y)),\atop
Ff(y)=-\frac{q^{5/2}}{1-q^2}(f(y)-f(q^2y))z^*.
\end{equation}
(\ref{action3}) allow one to "transfer" the structure of $U_q
\mathfrak{su}_{1,1}$-module from ${\rm Pol}({\mathbb C})_q$ onto $D({\mathbb
U})_q$.\medskip

{\bf Remark. } The functional (\ref{invint}) possesses the following
properties: for any $f\in D({\mathbb U})_q$, $\xi\in U_q\mathfrak{sl}_2$
\medskip

$1. \int \limits_{{\mathbb U}_q}f^*d\nu=\int \limits_{{\mathbb
U}_q}fd\nu$ (follows from the definition); \smallskip

$2. \int \limits_{{\mathbb U}_q}f^*fd\nu>0, f\ne0$ (\cite[Remark
3.6]{D2});\smallskip

$3. \int \limits_{{\mathbb U}_q}\xi fd\nu=\varepsilon(\xi)\cdot\int
\limits_{{\mathbb U}_q}fd\nu$ (\cite[Theorem 3.5]{D2}). \medskip

These properties allow one to regard the functional as a $q$-analogue of the
$SU(1,1)$-invariant integral.

Note (see \cite[Proposition 4.1]{D2}) that 1--3 imply: for any $f_1,f_2\in
D({\mathbb U})_q$, $\xi\in U_q \mathfrak{sl}_{2}$
\begin{equation}\label{unitary}
 (\xi f_1,f_2)=(f_1,\xi^* f_2),
\end{equation}
where $(f_1,f_2)\stackrel{\rm def}=\int \limits_{{\mathbb
U}_q}f_2^*f_1d\nu.$

The following formulae can be obtained
\begin{equation}\label{act1}
K^{\pm1}z^jf(y)=q^{\pm2j}z^jf(y),\quad
K^{\pm1}f(y)z^{*j}=q^{\mp2j}f(y)z^{*j},
\end{equation}
\begin{equation}\label{act2}
 E(z^jf(y))=-\frac{q^{1/2}}{1-q^2}z^{j+1}(f(y)-q^{2j}f(q^2y)),
\end{equation}
\begin{equation}\label{act3}
E(f(y)z^{*j})=-\frac{q^{1/2}}{1-q^2}((y-q^{-2j})f(y)+(1-y)f(q^{-2}y))
z^{*(j-1)},\quad j\geq1,
\end{equation}
\begin{equation}\label{act4}
F(z^{j}f(y))=-\frac{q^{5/2}}{1-q^2}z^{j-1}((y-q^{-2j})f(y)+(1-y)f(q^{-2}y))
,\quad j\geq1,
\end{equation}
\begin{equation}\label{act5}
 F(f(y)z^{*j})=-\frac{q^{5/2}}{1-q^2}(f(y)-q^{2j}f(q^2y))z^{*(j+1)}.
\end{equation}

Impose the notation $l_{i,j}$, \ $i=0,1,2, \ldots ,$ \
$j=0,\pm1,\pm2, \ldots , $ for the functional
$$\sum_{m>0}
z^m\psi_m(y)+\psi_0(y)+\sum_{m>0} \psi_{-m}(y)z^{*m} \mapsto
\psi_j(q^{2i})
$$
on the space $D({\mathbb U})'_q$ (see Section~1). Endow $D({\mathbb U})'_q$
with the weakest among the topologies in which all the linear functionals
$l_{i,j}$ are continuous. Obviously, $D({\mathbb U})_q$ is a dense subspace
in $D({\mathbb U})'_q$. As a straightforward consequence of
(\ref{act1})--(\ref{act5}) we get\medskip

\begin{proposition}\label{cont}
Any element $\xi\in U_q \mathfrak{sl}_{2}$ defines a continuous linear
operator $D({\mathbb U})_q\rightarrow D({\mathbb U})_q$ (here $D({\mathbb
U})_q$ is regarded as a topological vector space with the topology induced
by the topology on $D({\mathbb U})'_q$ described above).
\end{proposition}

\begin{corollary}\label{cont1} The $U_q \mathfrak{sl}_{2}$-action on
$D({\mathbb U})_q$ can be transferred by continuity onto the space
$D({\mathbb U})'_q$.
\end{corollary}

In fact the $U_q \mathfrak{sl}_{2}$-module $D({\mathbb U})'_q$ is a
covariant $D({\mathbb U})_q$-bimodule.

One can apply the above arguments to $D({\mathbb U})_q\otimes D({\mathbb
U})_q$, $D({\mathbb U}\times{\mathbb U})'_q$, $\{l_{i,j}\otimes l_{m,n}\}$
instead of $D({\mathbb U})_q$, $D({\mathbb U})'_q$, $\{l_{i,j}\}$ to make
$D({\mathbb U}\times{\mathbb U})'_q$ into a topological vector space and an
$U_q \mathfrak{sl}_{2}$-module. The continuity of the $U_q
\mathfrak{sl}_{2}$-action in $D({\mathbb U}\times{\mathbb U})'_q$ may be
proved just as in the case of $D({\mathbb U})'_q$.

The results listed below are proved in \cite{D2},\cite{D3}.
\medskip

\begin{proposition}  \label{invker}\cite[Proposition 4.5]{D2}. An
integral operator with a kernel $K$ is a morphism of the $U_q
\mathfrak{sl}_{2}$-module $D({\mathbb U})_q$ onto the $U_q
\mathfrak{sl}_{2}$-module $D({\mathbb U})'_q$ iff $K$ is an invariant.
\end{proposition}
\smallskip

\begin{proposition} \label{invgl}\cite[Section 6]{D3}. $G_l$ given
by (\ref{ker1}) is invariant for any $l\in{\mathbb C}$.
\end{proposition}\smallskip

\begin{proposition} \label{casimir}\cite[Theorem 4.3, Proposition
5.2]{D2} $\Delta_q$ being regarded as an operator $D({\mathbb
U})_q\rightarrow D({\mathbb U})'_q$ is a morphism of $U_q
\mathfrak{sl}_{2}$-modules.
 Moreover, $\Delta_q=q^{-1}\Omega$ where
\begin{equation}\label{casim}
 \Omega\stackrel{\rm
def}=FE+\frac{1}{(q^{-1}-q)^2}(q^{-1}K^{-1}+qK-q-q^{-1})
\end{equation}
is an element of the centre of $U_q \mathfrak{sl}_{2}$ called the Casimir
element. \end{proposition}\smallskip

\begin{proposition} \label{fzero}\cite[Theorem 3.9]{D2}.
 $f_0$ given by (\ref{f-0}) generates the $U_q
 \mathfrak{sl}_{2}$-module $D({\mathbb U})_q$.
\end{proposition}
\smallskip

\begin{corollary}\label{twomor} Let $A,B$ be morphisms of $U_q
\mathfrak{sl}_{2}$-modules $D({\mathbb U})_q\rightarrow D({\mathbb U})'_q$.
Then $A=B$ iff $Af_0=Bf_0$.
\end{corollary}

\section{\bf Proof of Theorem~2.2: reduction to the
 results of Section 3 about radial part of the quantum Laplacian}

Firstly it should be proved that
 the integral operators in the right-hand sides of (\ref{delta^-1})
 and (\ref{delta^-2}) are well defined (i.e., ${\mathbb G}_1$ and
${\mathbb G}_2$ do belong the space $D({\mathbb U}\times{\mathbb U})'_q$).
It could be done just as in the case of $G_l$ and $\hat{G}_N$ (see
Section~1) and we don't adduce such calculations.\medskip

\begin{lemma}\label{invgn} For any $N\in{\mathbb N}$ the kernel $\hat{G}_N$
given by (\ref{ker2}) is invariant.
\end{lemma}

{\bf Proof.} In our case the invariance of $\hat{G}_N$ means
$$E(\hat{G}_N)=F(\hat{G}_N)=(K^{\pm 1}-1)(\hat{G}_N)$$ and follows
from the continuity of the $U_q \mathfrak{sl}_{2}$-action in $D({\mathbb
U}\times{\mathbb U})'_q$ and Proposition \ref{invgl}. \hfill $\square$

\begin{lemma}\label{gf0}
\begin{equation}\label{g1f0}
 \int \limits_{{\mathbb U}_q}{\mathbb
G}_1f_0d\nu=g_1(y),\quad \int \limits_{{\mathbb U}_q}{\mathbb
G}_2f_0d\nu=g_2(y),
\end{equation}
where $g_1(y)$, $g_2(y)$ are given by (\ref{g1}) and (\ref{g2})
respectively.
\end{lemma}\smallskip

{\bf Proof.} (\ref{g1f0}) reduce to (\ref{ker1}), (\ref{ker3}) and the
following equalities:\medskip

1. $z^*\cdot f_0=0\quad $ (see \cite[Proposition 3.1]{D2});\smallskip

2. $f(y)\cdot f_0=f(1)\quad $ (follows from the definition (\ref{f-0}) of
$f_0$); \smallskip

3. $\int \limits_{{\mathbb U}_q}z^kf_0d\nu=
\begin{cases}
1-q^2, & k=0,\\ 0, & k=1,2,\ldots
\end{cases} $
(follows from the definition (\ref{invint}) of the integral). \hfill
$\square$
\smallskip

Thus we have proved (Proposition \ref{invgl}, Lemma \ref{invgn}) that the
operators in the right-hand sides of (\ref{delta^-1}) and (\ref{delta^-2})
are morphisms of $U_q \mathfrak{sl}_{2}$-modules. By Lemma \ref{gf0} and
Theorem \ref{g1,g2}
$$\int \limits_{{\mathbb
U}_q}{\mathbb G}_1f_0d\nu=\Delta_q^{-1}f_0,
$$
$$\int \limits_{{\mathbb U}_q}{\mathbb
G}_2f_0d\nu=\Delta_q^{-2}f_0.$$

By Corollary \ref{twomor} to complete the proof of Theorem
\ref{princresult} it suffices now to prove the following lemma.\medskip

\begin{lemma}\label{morph} $\Delta_q^{-1}$ and $\Delta_q^{-2}$ being
regarded as operators from $D({\mathbb U})_q$ onto $D({\mathbb U})'_q$ are
morphisms of $U_q \mathfrak{sl}_{2}$-modules.
\end{lemma}

{\bf Proof.} Let $t_\phi$ be the automorphism of the algebra ${\rm
Pol}({\mathbb C})_q$ given by
$$ t_\phi(z)=e^{i\phi}z,\quad t_\phi(z^*)=e^{-i\phi}z^*.
$$
Impose the same notation $t_\phi$ for the automorphism
$$\sum_{m>0}
z^m\psi_m(y)+\psi_0(y)+\sum_{m>0} \psi_{-m}(y)z^{*m}
$$
$$
 \mapsto\sum_{m>0}e^{im\phi} z^m\psi_m(y)+\psi_0(y)+\sum_{m>0}e^{-im\phi}
\psi_{-m} (y)z^{*m}$$ of the algebra $D({\mathbb U})_q$.

Obviously, each operator $t_\phi$ can be extend to a unitary operator
$L^2(d\nu)_q \rightarrow L^2(d\nu)_q$ and thus we obtain a unitary
representation of the group $U(1)$.

Let's prove that for any $\phi$
$$t_\phi\cdot\Delta_q=\Delta_q\cdot t_\phi.$$
Indeed, let $L_n$ be the subspace in $D({\mathbb U})_q$ of function of the
form $z^nf(y)$ (for $n>0$), $f(y)$ (for $n=0$) or $f(y)z^{*n}$ (for $n<0$).
Then it is obvious that
$$L_n=\{f\in D({\mathbb U})_q:t_\phi(f)=e^{in\phi}f, \phi\in[0;2\pi)\}=
\{f\in D({\mathbb U})_q:K(f)=q^{2n}f\}.
$$
 The latter equality and Proposition \ref{casimir}
imply $\Delta_q(L_n)\subset L_n$ and thus
$$t_\phi(\Delta_q(f))=\Delta_q(t_\phi(f))$$ for any $t_\phi$ and $f\in
D({\mathbb U})_q$.

Let $\overline{L_n}$ be the closure of $L_n$ in $L^2(d\nu)_q$. It is evident
that $\Delta_q(\overline{L_n})\subset \overline{L_n}$ and, moreover,
$\Delta_q^{-1}(\overline{L_n})\subset \overline{L_n}$ (this follows from the
invertibility of $\Delta_q$).

Denote by $L^2(d\nu)^{fin}_q$ the space $\bigoplus_{n\in {\mathbb
Z}}\overline{L_n}$ (note that $\overline{L_n}\perp\overline{L_m}$ for $n\ne
m$). We have established that
$$\Delta_q(L^2(d\nu)^{fin}_q)\subset L^2(d\nu)^{fin}_q,$$ and
$$\Delta_q^{-1}(L^2(d\nu)^{fin}_q)\subset L^2(d\nu)^{fin}_q.$$

Obviously, $D({\mathbb U})_q\subset L^2(d\nu)^{fin}_q$ and therefore
$$\Delta_q^{-m}(D({\mathbb U})_q)\subset L^2(d\nu)^{fin}_q$$
for any $m\in{\mathbb N}$.

To complete the proof of Lemma \ref{morph} it suffices to prove that
$L^2(d\nu)^{fin}_q$ is an $U_q{\mathfrak sl}_2$-submodule in $D({\mathbb
U})'_q$. Thus we have to verify inclusions
\begin{equation}\label{inc1}
 K(L^2(d\nu)^{fin}_q)\subset L^2(d\nu)^{fin}_q,
\end{equation}
\begin{equation}\label{inc2}
 E(L^2(d\nu)^{fin}_q)\subset L^2(d\nu)^{fin}_q,
\end{equation}
\begin{equation}\label{inc3}
 F(L^2(d\nu)^{fin}_q)\subset L^2(d\nu)^{fin}_q.
\end{equation}
(\ref{inc1}) is evident. Let us prove (\ref{inc2}) ((\ref{inc3}) can be
proved in a similar way).

Formulae (\ref{act3}), (\ref{act5}), imply $E(L_n)\subset L_{n+1}$ and we
need to prove that $E$ is extendable onto $\overline{L_n}$. Let
$f\in\overline{L_n}$. Then (see (\ref{unitary}))
$$ (Ef,Ef)=(f,E^*Ef)\stackrel{({\rm see}(\ref{invol}))}=-(f,KFEf)$$
$$
\stackrel{({\rm
see}(\ref{casim}))}=-(f,K\cdot(\Omega-\frac{q^{-1}K^{-1}+
qK-q^{-1}-q}{(q^{-1}-q)^2})f)$$
$$=-(f,K\Omega f)+(f,\frac{q^{-1}+
qK^2-(q^{-1}+q)K}{(q^{-1}-q)^2}f)$$
$$=-q^{2n}(f,\Omega
f)+\frac{q^{-1}+ q^{4n+1}-(q^{-1}+q)q^{2n}}{(q^{-1}-q)^2}(f,f)$$
$$
\stackrel{({\rm Proposition~4.5})}=-q^{2n+1}(f,\Delta_q
f)+ \frac{q^{-1}+ q^{4n+1}-(q^{-1}+q)q^{2n}}{(q^{-1}-q)^2} (f,f).$$

So the boundedness of $\Delta_q$ allows one to establish the boundedness of
$E:L_n \rightarrow L_{n+1}$. This completes the proof of Lemma \ref{morph}
and thus of Theorem \ref{princresult}. \medskip

{\bf Acknowledgement.} I am very grateful to L.L.~Vaksman who explained me
general ideas \cite{V}, \cite{V1} of producing invariant kernels for quantum
homogeneous spaces.
\end{document}